\documentclass[12pt]{article}
\usepackage{amsfonts,amsmath}
\usepackage{tikz}
\usepackage{pgfplots}
\usepackage{pgfplotstable}
\usepackage{booktabs}
\usepackage{float}
\usepackage{siunitx}

\usepackage{a4wide,epsfig}
\usepackage[utf8]{inputenc}   % damit man Umlaute direkt eingeben kann und diese erkannt werden.

\usepackage{todonotes}
\usepackage{algorithm}
\usepackage{algpseudocode}
\usepackage{pgfplots}
\usepackage{pgfplotstable}
\usepackage{booktabs}

\numberwithin{equation}{section}

\begin{document}

\setcounter{page}{1}

\title{Efficient Solution of State-Constrained Distributed Parabolic Optimal
Control Problems}
\author{Richard~L\"oscher, Michael~Reichelt, Olaf~Steinbach}
\date{Institut f\"ur Angewandte Mathematik, TU Graz, \\ 
Steyrergasse 30, 8010 Graz, Austria}

\maketitle

\begin{abstract}
  We consider a space-time finite element method for the
  numerical solution of a distributed tracking-type optimal control problem
  subject to the heat equation with state constraints. The cost or
  regularization term is formulated in an anisotropic Sobolev norm for the
  state, and the optimal state is then characterized as the unique solution
  of a first kind variational inequality. We discuss an efficient
  realization of the anisotropic Sobolev norm in the case of a space-time
  tensor-product finite element mesh, and the iterative solution of the
  resulting discrete variational inequality by means of a semi-smooth
  Newton method, i.e., using an active set strategy.
\end{abstract}

\section{Introduction}
Optimal control problems arise naturally in a wide range of applications,
e.g., \cite{Reichelt:HerzogHeinkenschlossKaliseStadlerTrelat:2022}.
In this work, we consider tracking type optimal control problems to minimze 
\begin{equation}\label{Reichelt:functional}
  {\mathcal{J}}(u_\varrho,z_\varrho) =
  \frac{1}{2} \int_0^T \int_\Omega [u_\varrho(x,t)-\overline{u}(x,t)]^2 \,
  dx \, dt + \frac{1}{2} \, \varrho \, \| z_\varrho \|^2_Z
\end{equation}
subject to the Dirichlet boundary value problem for the heat equation,
\begin{equation}\label{Reichelt:heat equation}
  \begin{array}{rclcl}
    \partial_t u_\varrho(x,t) - \Delta_x u_\varrho(x,t)
    & = & z_\varrho (x,t) && \mbox{for} \; (x,t) \in
                             Q := \Omega \times (0,T),\\[1mm]
    u_\varrho(x,t) & = & 0 && \mbox{for} \; (x,t) \in \Sigma :=
                              \partial\Omega \times (0,T), \\[1mm]
    u_\varrho(x,0) & = & 0 && \mbox{for} \; x \in \Omega,
  \end{array}
\end{equation}
where $\Omega \subset \mathbb{R}^n, n=2,3$ is a bounded Lipschitz domain,
and $T>0$ is a finite time horizon. In \eqref{Reichelt:functional}, we aim to
approximate a given, probably discontinuous, target
$\overline{u} \in L^2(Q)$ by a more regular function $u_\varrho \in X$
satisfying the heat equation \eqref{Reichelt:heat equation} with the control
$z_\varrho$ as right hand side. For $z_\varrho \in Z$ we write
\eqref{Reichelt:heat equation} as operator equation to find
$u_\varrho \in X$ such that $B u_\varrho=z_\varrho$ is satisfied, and
we assume that $B : X \to Z$ defines an isomorphism.
Instead of \eqref{Reichelt:functional} we then consider the reduced
functional for $u_\varrho \in X$,
\begin{equation}\label{Reichelt:reduced functional}
  \widetilde{\mathcal{J}}(u_\varrho) =
  \frac{1}{2} \int_0^T \int_\Omega [u_\varrho(x,t)-\overline{u}(x,t)]^2 \,
  dx \, dt + \frac{1}{2} \, \varrho \, \| B u_\varrho \|^2_Z \, .
\end{equation}
The most standard choice is $Z=L^2(Q)$, and hence,
using $Y = L^2(0,T;H^1_0(\Omega))$,
\[
  X := \{ u \in Y :
    \partial_t u \in Y^*, u(x,0) = 0, \; x \in \Omega,
    \partial_t u - \Delta_x u \in L^2(Q) \} .
\]
For a non-conforming discretization of the resulting reduced optimality system,
using piecewise linear continuous space-time finite element functions
defined with respect to a simplicial decomposition of the space-time
domain $Q$, see \cite{Reichelt:LangerSteinbachTroeltzschYang:2021a}.
This finite element discretization becomes conforming if we consider
$Z = Y^* = L^2(0,T;H^{-1}(\Omega))$ and
$X := \{ u \in Y : \; \partial_t u \in Y^*, u(x,0) = 0, \; x \in \Omega \}$,
see \cite{Reichelt:LangerSteinbachTroeltzschYang:2021b}. Moreover, following
\cite{Reichelt:SteinbachZank:2020}, we can consider the Dirichlet
boundary value problem \eqref{Reichelt:heat equation} in anisotropic
Sobolev spaces, i.e., $X = H^{1,1/2}_{0;0,}(Q)$, and
$Z = [H^{1,1/2}_{0;,0}(Q)]^*$, see also
\cite{Reichelt:löscher2024optimalcomplexitysolutionspacetime}.
In all of these cases, the space-time finite element discretization
of the reduced optimality system results, after eliminating the discrete
adjoint state, in algebraic systems of linear equations of the form
$[M_h + \varrho B_h^\top A_h^{-1} B_h] \underline{u} = \underline{f}$.
Here, $B_h$ is the space-time finite element matrix which is related
to the Dirichlet boundary value problem \eqref{Reichelt:heat equation},
$M_h$ is the mass matrix coming from the first part in
\eqref{Reichelt:reduced functional}, and $A_h$ is a space-time finite
element matrix in order to realize a norm in $Z^*$. In any case, the
structure of the Schur complement matrix $S_h := B_h^\top A_h^{-1} B_h$ may
complicate an efficient iterative solution of the discrete reduced
optimality system, due to the involved inversion of $A_h$. Instead of the
Schur complement system we may also solve an equivalent block skew-symmetric
but positive definite system, see, e.g.
\cite{Reichelt:LangerLoescherSteinbachYang:2023} in the case of the
Poisson equation with $L^2$ regularization.

To avoid the application
of the Schur complement $S_h$ we may replace $S_h$ by any spectrally
equivalent space-time finite element stiffness matrix $D_h$ which
realizes a norm in $X = H^{1,1/2}_{0;0,}(Q)$. Instead of
\eqref{Reichelt:reduced functional} we therefore consider the
minimization of 
\begin{equation}\label{Reichelt:state functional}
  \widehat{\mathcal{J}}(u_\varrho) =
  \frac{1}{2} \, \| u_\varrho - \overline{u} \|_{L^2(Q)}^2 +
  \frac{1}{2} \, \varrho \, \| u_\varrho \|_{H^{1,1/2}_{0;0,}(Q)}^2
\end{equation}
which results in the determination of the optimal state $u_\varrho$
from which we can compute the optimal control $z_\varrho = B u_\varrho$
by some post processing. While in
\cite{Reichelt:löscher2024optimalcomplexitysolutionspacetime} we have
considered the minimization of \eqref{Reichelt:state functional}
without additional state or control constraints, here we include state
constraints, see also
\cite{Reichelt:gangl2023regularizationfiniteelementerror}, i.e.,
\begin{equation}\label{Reichelt:state constraints}
  u_\varrho \in \mathcal{K} = \lbrace u \in X \mid u_- \leq u \leq u_+ \text{
    a.e. in } Q\rbrace \, ,
\end{equation}
where $u_\pm \in X \cap C(\overline{Q})$ are given continuous
barrier functions, and we assume $0 \in {\mathcal{K}}$.

The norm in $X=H^{1,1/2}_{0;0,}(Q)$ is realized via 
\begin{equation}\label{Reichelt:Def Norm}
  \| u \|_{H^{1,1/2}_{0;0,}(Q)}^2 :=
    \langle \partial_t u , {\mathcal{H}}_T u \rangle_Q +
    \| \nabla_x u \|^2_{L^2(Q)} =: \langle D u , u \rangle_Q \, ,
\end{equation}
where $\mathcal{H}_T$ is the modified Hilbert transform
\cite{Reichelt:SteinbachZank:2020}, that only acts in the temporal direction, 
\[
  {\mathcal{H}}_Tu(x,t) = \sum\limits_{k=0}^\infty u_k(x) \cos \left( \Big(
    \frac{\pi}{2} + k \pi \Big) \frac{t}{T} \right),
\]
and the Fourier coefficients are given as
\[
  u_k(x) = \frac{2}{T} \int_0^T u(x,t)  \sin \left( \Big(
    \frac{\pi}{2} + k \pi \Big) \frac{t}{T} \right) \, dt .
\]
Note that we can write, using $B u_\varrho = z_\varrho$,
\[
  \| u_\varrho \|^2_{H^{1,1/2}_{0;0,}(Q)} =
  \langle D u_\varrho , u_\varrho \rangle_Q =
  \| u_\varrho \|_D^2 = \| B^{-1} z_\varrho \|^2_D,
\]
which defines an equivalent norm for the control $z_\varrho$
in $[H^{1,1/2}_{0;,0}(Q)]^*$.

While in the unconstrained case the minimization of
\eqref{Reichelt:state functional} results in the gradient
equation $u_\varrho + \varrho D u_\varrho = \overline{u}$, in the case of
state constraints we find the optimal state
$u_\varrho \in {\mathcal{K}}$ as unique solution of the first kind
variational inequality
\begin{equation}\label{Reichelt:variational inequality}
  \langle u_\varrho , v - u_\varrho \rangle_{L^2(Q)} +
  \varrho \, \langle D u_\varrho , v - u_\varrho \rangle_Q \geq
  \langle \overline{u} , v - u_\varrho \rangle_{L^2(Q)}
  \quad \mbox{for all} \; v \in {\mathcal{K}}.
\end{equation}
For the numerical solution of \eqref{Reichelt:variational inequality}
we define a space-time tensor product
finite element space $X_h := W_{h_x} \otimes V_{h_t} \subset X$, where
$W_{h_x} = \mbox{span} \{ \psi_i \}_{i=1}^{M_x} \subset H^1_0(\Omega)$
is the spatial finite element space of piecewise linear basis functions
$\psi_i$ which are defined with respect to some admissible and globally
quasi-uniform finite element mesh with spatial mesh size $h_x$. Further,
$V_{h_t} := S_{h_t}^1(0,T) \cap H^{1/2}_{0,}(0,T) =
\mbox{span} \{ \varphi_k \}_{k=1}^{N_t}$ is the space of piecewise linear
functions, which are defined with respect to a uniform finite element
mesh with temporal mesh size $h_t$. With this we define
\[
  {\mathcal{K}}_h := \lbrace u_h \in X_h : I_h u_- \leq u_h \leq I_h u_+
  \; \mbox{in} \; Q \rbrace,
\]
and we consider the Galerkin discretization of the variational inequality
\eqref{Reichelt:variational inequality} to find
$u_{\varrho h} \in {\mathcal{K}}_h$ such that
\begin{equation}\label{Reichelt:variational inequality FEM}
  \langle u_{\varrho h} , v_h - u_{\varrho h} \rangle_{L^2(Q)} +
  \varrho \, \langle D u_{\varrho h} , v_h - u_{\varrho h} \rangle_Q \geq
  \langle \overline{u} , v_h - u_{\varrho h} \rangle_{L^2(Q)}, \;
  \forall v_h \in {\mathcal{K}}_h.
\end{equation}
When assuming $\overline{u} \in {\mathcal{K}} \cap H^{2,1}(Q)$ and
$u_\pm \in H^{2,1}(Q)$ we can prove the following error estimate,
see \cite{Reichelt:gangl2023regularizationfiniteelementerror,
  Reichelt:löscher2024optimalcomplexitysolutionspacetime},
\begin{eqnarray*}
  && \| u_{\varrho h} - \overline{u} \|^2_{L^2(Q)} \\
  && \leq \,
  c \, \Big[ h_t^2 + h_x^4 + \varrho (h_t + h_x^2) + \varrho^2 \Big] \, \Big[ \| \overline{u} \|^2_{H^{2,1}(Q)} +
  \| u_+ \|^2_{H^{2,1}(Q)} + \| u_- \|^2_{H^{2,1}(Q)} \Big] .
\end{eqnarray*}
This error estimate motivates the particular choices
$\varrho = h_x^2$ and $h_t = h_x^2$ in order to conclude
\[
\| u_{\varrho h} - \overline{u} \|^2_{L^2(Q)} \, \leq \,
  c \, h_x^4\, \Big[ \| \overline{u} \|^2_{H^{2,1}(Q)} +
  \| u_+ \|^2_{H^{2,1}(Q)} + \| u_- \|^2_{H^{2,1}(Q)} \Big] .
\]
Note that for a more regular target, e.g., for
$\overline{u} \in {\mathcal{K}} \cap H^2(Q)$, we can also use $h_t=h_x$
to obtain the same estimate.

In this note we aim for an efficient iterative solution of the discrete
variational inequality \eqref{Reichelt:variational inequality FEM}.
Albeit $D$ as defined in \eqref{Reichelt:Def Norm} is a non-local operator,
since we are using a space-time tensor product ansatz space,
we will demonstrate in Section \ref{sec:discretization} an efficient
matrix free discretization of \eqref{Reichelt:variational inequality FEM}.
These results are then used for an iterative solution of the 
discretized optimality system. This has the advantage, that we can augment
this operator to yield an efficient algorithm for the active set strategy
later on. The description of the active set strategy with the augmented
operator is presented in Section \ref{sec:active set}. Numerical results
are presented in Section \ref{sec:numerical-results}.

\section{Discretization}
\label{sec:discretization}
The discrete variational inequality \eqref{Reichelt:variational inequality FEM}
can be written in the following form to find
$\underline{u} \in {\mathbb{R}}^{N_t \cdot M_x}
\leftrightarrow u_h \in {\mathcal{K}}_h$ such that
\begin{equation}\label{Reichelt:discrete VI}
  (K_h \underline{u} , \underline{v} - \underline{u} ) \geq
  (\underline{f} , \underline{v} - \underline{u}) \quad
  \mbox{for all} \; \underline{v} \in {\mathbb{R}}^{N_t \cdot M_x}
  \leftrightarrow v_h \in {\mathcal{K}}_h,
\end{equation}
where $K_h$ is the space-time finite element approximation of
$I+\varrho D$, and $\underline{f}$ is the load vector 
which is related to the given target $\overline{u}$.
Since we are using a space-time tensor product ansatz space
$X_h$, we can write the space-time finite element matrix as
\begin{equation}\label{Reichelt:stiffness matrix}
  K_h = M_{h_t} \otimes M_{h_x} + \varrho \, \Big[
  A_{h_t} \otimes M_{h_x} + M_{h_t} \otimes A_{h_x} \Big] \in
  {\mathbb{R}}^{N_t \cdot M_x \times N_t \cdot M_x} ,
\end{equation}
where
\[
  A_{h_t}[j,i] = \langle \partial_t \varphi_i , {\mathcal{H}}_T
  \varphi_j \rangle_{L^2(0,T)}, \;
  M_{h_t}[j,i] = \langle \varphi_i , \varphi_j \rangle_{L^2(0,T)},
  \; i,j=1,\ldots,N_t,
\]
\[
  A_{h_x}[\ell,k] = \langle \nabla_x \psi_k , \nabla_x \psi_\ell
  \rangle_{L^2(\Omega)}, \;
  M_{h_x}[\ell,k] = \langle \psi_k , \psi_\ell \rangle_{L^2(\Omega)},
  \; k,\ell=1,\ldots,M_x.
\]
For the optimal choice $\varrho = h_x^2$, the matrix $K_h$ is spectrally
equivalent to the space-time mass matrix $M_{h_t} \otimes M_{h_x}$, allowing
for a simple diagonal preconditioning when inverting $K_h$ by applying
a preconditioned conjugate gradient scheme. In order to realize the
matrix-vector product $\underline{w} = K_h \underline{v}$ efficiently,
in \cite{Reichelt:löscher2024optimalcomplexitysolutionspacetime} we
have considered the generalized eigenvalue problem 
\begin{equation}\label{gen EVP}
  A_{h_t} \underline{v} = \lambda \, M_{h_t} \underline{v} .
\end{equation}
With this we are able to transform any vector
$\underline{v} \in {\mathbb{R}}^{N_t}$ into the eigenvector basis
of $\left(A_{h_t}, M_{h_t} \right)$ with $\mathcal{O}(N_t \log N_t)$ effort,
which translates to $\mathcal{O}(N_t M_x \log N_t)$ for any vector in
${\mathbb{R}}^{N_t \cdot M_x}$. The temporal transformation matrix into the
eigenvector basis is denoted by $C_{h_t}^{-1}$. The respective transformation
on ${\mathbb{R}}^{N_t \cdot M_X}$ is then given by $C_{h_t}^{-1} \otimes I_{h_x}$. 
For the application $\underline{w} = K_h \underline{v}$ we then conclude
$\underline{w} = (M_{h_t}C_{h_t}\otimes I_x) \hat{K}_h
(C_{h_t}^{-1} \otimes I_x) \underline{v}$,
where $\hat{K}_h$ is the representation of $K_h$ in the temporal eigenbasis, i.e., 
\[
  \hat{K}_h = (C_{h_t}^{-1} M_{h_t}^{-1} \otimes I_{h_x}) K_h
  (C_{h_t}\otimes I_{h_x}) 
  =
  I_{h_t} \otimes M_{h_x} + \varrho
  (\Lambda_{h_t} \otimes I_{h_x} + I_{h_t} \otimes A_{h_x}),
\]  
and $\Lambda_{h_t}$ is the diagonal matrix of the generalized eigenvalues
$\lambda_i$ of \eqref{gen EVP}. The application of $\hat{K}_h$ has effort
$\mathcal{O}(N_t M_x)$, so we end up with the overall effort of
$\mathcal{O}(N_t M_x \log N_t)$ for the matrix free application of $K_h$.
This is a significant improvement over the $\mathcal{O}(N_t^2 M_x)$ effort
for the direct application of $K_h$. Together with the feasibility of the
conjugate gradient method with diagonal preconditioning, this yields a
quasi optimal solver. Additionally, shared memory parallelization can be
easily implemented, due to the Kronecker product structure of the arising
matrices.

\section{Semi-smooth Newton method}\label{sec:active set}
For the solution $\underline{u}$ of the discrete variational inequality
\eqref{Reichelt:discrete VI} we define $\underline{\lambda} :=
K_h \underline{u} - \underline{f}$. By 
${\mathcal{I}}_{A,\pm}$ we denote the index set of all active nodes
where $u_j = u_{\pm,j} := u_\pm(x_j,t_j)$, while the complementary set is
called the inactive set ${\mathcal{I}}_{I}$.
Then there hold the discrete complementarity conditions
\[
  \lambda_j = 0, \; u_{-,j} < u_j < u_{+,j} \; \mbox{for} \;
  j \in {\mathcal{I}}_I,
  \;
  \lambda_j \leq 0 \; \mbox{for} \; j \in {\mathcal{I}}_{A,+}, \;
  \lambda_j \geq 0 \; \mbox{for} \; j \in {\mathcal{I}}_{A,-},
\]
which are equivalent to
\[
  \lambda_j = \min \{ 0 , \lambda_j +c(u_{+,j}-u_j) \} +
  \max \{ 0 , \lambda_j + c (u_{-,j}-u_j) \} , \quad c > 0.
\]
Hence we have to solve a system
$\underline{F}(\underline{u},\underline{\lambda})= \underline{0}$
of (non)linear equations
\begin{eqnarray*}
  \underline{F}_1(\underline{u},\underline{\lambda})
  & = &
  K_h \underline{u} - \underline{\lambda} - \underline{f} =
  \underline{0}, \\
  \underline{F}_2(\underline{u},\underline{\lambda})
  & = &
  \underline{\lambda} - \min \{0,\underline{\lambda}+c(\underline{u}_+-
  \underline{u})) \} - \max \{0,\underline{\lambda}+c(\underline{u}_--
  \underline{u})) \} = \underline{0} .
\end{eqnarray*}
One step of the semi-smooth Newton method
\cite{Reichelt:Kunisch:2002}
for the iterative solution of this system reads
\begin{align}
  \begin{pmatrix}
    \underline{u}^{k+1} \\%
    \underline{\lambda}^{k+1}
  \end{pmatrix} = %
  \begin{pmatrix}
    \underline{u}^{k} \\%
    \underline{\lambda}^{k}
  \end{pmatrix}
  - \left[D\underline{F}(\underline{u}^k, \underline{\lambda}^k)\right]^{-1}
  \underline{F}(\underline{u}^k, \underline{\lambda}^k),
  \label{eq:semismooth:newton}
\end{align}
where $D\underline{F}$ is the Jacobian of $\underline{F}$ in the sense of
slant derivatives. For any iterate $(\underline{u}^k,\underline{\lambda}^k)$
we define the active set ${\mathcal{I}}_A^k$ as well as the inactive set
${\mathcal{I}}_I^k$. From the second line in \eqref{eq:semismooth:newton}
we first conclude
\begin{align}
  \begin{cases}
    \lambda^k_j + c (u_{-,j} - u^k_j ) > 0, \\%
    \lambda^k_j + c (u_{+,j} - u^k_j ) < 0, \\%
    \text{else},
  \end{cases}
  \Longrightarrow
  \begin{cases}
    u_j^{k+1} := u_{-,j}, \\%
    u_j^{k+1} := u_{+,j}, \\%
    \lambda_j^{k+1} := 0.
  \end{cases} \label{eq:active:set}
\end{align}
For the remainder, we first introduce a split into active and inactive
parts of all quantities,
\begin{align}
  \underline{u}^{k+1} = \underline{u}^{k+1}_A \oplus \underline{u}^{k+1}_I,
  \quad
  \underline{\lambda}^{k+1} = -\underline{\lambda}^{k+1}_A \oplus
  \underline{\lambda}^{k+1}_I,
\end{align}
where $\lambda^{k+1}_{I,j} = 0$ for $j \in \mathcal{I}^{k+1}_I$ and
$u^{k+1}_{A,j}$ is set according to \eqref{eq:active:set} for
$j \in \mathcal{I}^{k+1}_A$. Substituting this into the remaining first
equation of \eqref{eq:semismooth:newton} yields
\begin{align}
  \widetilde{K}_h^{k+1}
  \begin{pmatrix}
    \underline{u}^{k+1}_I \\[1mm]
    \underline{\lambda}^{k+1}_A
  \end{pmatrix}
  :=  K_h \underline{u}^{k+1}_I + \underline{\lambda}^{k+1}_A
  = \overline{\underline{u}} - K_h \underline{u}^{k+1}_A
  =: \underline{f}^{k+1}. \label{eq:active:set:linear:problem}
\end{align}
The left-hand side is a linear map
$\widetilde{K}_h^{k+1} : \mathbb{R}^{N_t M_x} \to \mathbb{R}^{N_t M_x}$ for
any choice of an active set. Due to the orthogonality of the splitting,
$K_h \underline{u}^{k+1}_I + \underline{\lambda}^{k+1}_A$ is also an
orthogonal sum and hence $\widetilde{K}_h^{k+1}$ inherits positive
definiteness from $K_h$ and the identity. After a correct reordering
of the indices, the matrix is even block diagonal,
\begin{align}
  \widetilde{K}_h^{k+1} = \begin{pmatrix}
    R_I^{k+1} K_h^{k+1} P_I^{k+1} & 0 \\[1mm]
    0 & R_A^{k+1} I_h P_A^{k+1}
  \end{pmatrix},
\end{align}
where $P^{k+1}_{I/A}$ and $R^{k+1}_{I/A}$ are the canonical prolongation
(by zero) and restriction with respect to their subscript index set.
Now feasibility of diagonal preconditioning is evident from its feasibility
for $K_h$. In practice the splitting is performed by setting vanishing
components to zero. Then all matrices can be applied as in the unconstrained
case. As $\widetilde{K}_h$ is not available as a matrix, we take the diagonal
of the spectrally equivalent space-time mass matrix for preconditioning,
where we set entries corresponding to the active set to one. The semi-smooth
Newton method terminates, when the index sets are no longer changing.
The algorithm for the overall procedure is summarized in
Algorithm \ref{alg:active:set}. Note, that the simple stopping criterion is
no longer applicable, when, e.g., underrelaxation or line-search is used.
Then one needs to add additional stopping criteria. Underrelaxation is
discussed in the next chapter.

\begin{algorithm}
  \caption{Active set strategy in the case of state constraints}
  \label{alg:active:set}
  \begin{algorithmic}
  \State choose initial guesses for $\underline{u}^0$ and $\underline{\lambda}^0$
  \For {$k = 0,1,2,\ldots$}
    \State compute $\mathcal{I}_A^{k+1}$ and $\mathcal{I}_I^{k+1}$ according to \eqref{eq:active:set}
    \If{$k>0$ and $\mathcal{I}_A^{k+1} = \mathcal{I}_A^{k}$} \\
      \Return $\left(\underline{u}^{k}, \underline{\lambda}^k\right)$ 
    \EndIf
    \State solve \eqref{eq:active:set:linear:problem} for $\underline{u}^{k+1}_I$ and $\underline{\lambda}^{k+1}_A$ using matrix free CG with diagonal preconditioning
  \EndFor
  \end{algorithmic}
\end{algorithm}

\section{Numerical Results}\label{sec:numerical-results}
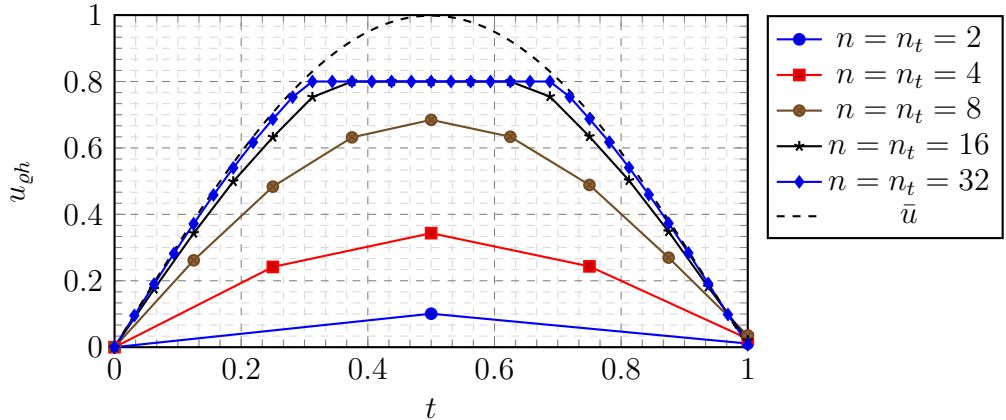
\begin{figure}[ht]
  \centering
  \begin{tikzpicture}
      \begin{axis}[
          xlabel={$t$},
          ylabel={$u_{\varrho h}$},
          grid=both,
          width=10cm,
          height=6cm,
          grid style={dashed, gray!30},
          major grid style={black!50},
          minor tick num=5,
          thick,
          cycle list name=auto,
          ymin=0, ymax=1, % Adjust these limits as needed
          xmin=0, xmax=1.0,
          legend pos=outer north east,
      ]
      \addplot %
      table[x=t,y=u,col sep=comma] {simulation_data/trajectories/trajectory_constrained_3d_refinement_0.csv};
      \addlegendentry{$n=n_t=2$}
      
      \addplot %
      table[x=t,y=u,col sep=comma] {simulation_data/trajectories/trajectory_constrained_3d_refinement_1.csv};
      \addlegendentry{$n=n_t=4$}
      
      \addplot %
      table[x=t,y=u,col sep=comma] {simulation_data/trajectories/trajectory_constrained_3d_refinement_2.csv};
      \addlegendentry{$n=n_t=8$}
      
      \addplot %
      table[x=t,y=u,col sep=comma] {simulation_data/trajectories/trajectory_constrained_3d_refinement_3.csv};
      \addlegendentry{$n=n_t=16$}
      
      \addplot %
      table[x=t,y=u,col sep=comma] {simulation_data/trajectories/trajectory_constrained_3d_refinement_4.csv};
      \addlegendentry{$n=n_t=32$}

      % Add a plot of u(x=0.51, y=0.51, z=0.51, t) 
      \addplot[domain=0:1, samples=200, thick, dashed] {sin(deg(pi*0.51))*sin(deg(pi*0.51))*sin(deg(pi*0.51))*sin(deg(pi*x))};
      \addlegendentry{$\bar{u}$}

      \end{axis}
  \end{tikzpicture}
  \caption{Plot of the constrained solution $u_{\varrho h}$ along the line $x_1=x_2=x_3=0.51$. }\label{fig:constrained_solution}
\end{figure}

As an illustrative example, we consider $\Omega = (0,1)^3$ and $T=1$, i.e.,
the space-time cylinder $Q = \Omega \times (0,T)= (0,1)^4$. We decompose
$\Omega$ into a shape regular and simplicial globally quasi-uniform
mesh $\Omega_h$ with $n_x$ elements
of mesh size $h_x$, and $(0,T)$ into a regular mesh $\mathcal{T}_h$
consisting of $n_t$ equidistant intervals of length $h_t=1/n_t$. Then we
define the discrete space $X_h = W_{h_x} \otimes V_{h_t}$ with
$W_{h_x} = S^1_{h_x}(\Omega_h) \cap H^1_0(\Omega)$ and $V_{h_t} = S^1_{h_t}(\mathcal{T}_h)
\cap H^{1/2}_{0,}(0,T)$, where $S^1_{h}$ denotes the space of piecewise linear
continuous functions on the respective mesh.
As target function we consider
\begin{align}
  \overline{u}(x,t) =
  \sin \left(\pi x_1\right) \sin \left(\pi x_2\right)
  \sin \left(\pi x_3\right) \sin \left(\pi t\right) \in C^\infty(Q) \cap X.
\end{align}
As we have extensively validated the optimality of the strongly related
solver for the unconstrained case in
\cite{Reichelt:löscher2024optimalcomplexitysolutionspacetime}, here we
will only present the more interesting case with state constraints, choosing
$u_- \equiv 0$ and $u_+ \equiv 0.8$. This implies that
$\overline{u} \notin {\mathcal{K}}$. But all the given error estimates
remain valid when replacing $\overline{u}$ by its projection
$P_{\mathcal{K}}\overline{u} \in {\mathcal{K}}$. As initial guess in the
semi-smooth Newton method we consider
$\underline{u}^0 = \frac{1}{2} (\underline{u}_- + \underline{u}_+)$
and $\underline{\lambda}^0 = K_h \underline{u}^0 - \underline{f}$.
Further we apply underrelaxation with $\omega_N = 0.1$ in the semi-smooth
Newton method with the parameter $c=1$. We stop the non-linear solver if
the active sets do not change anymore and subsequent iterates satisfy
\begin{align*}
  \lVert \underline{u}^{k+1}-\underline{u}^k \rVert_{\infty} +
  \lVert \underline{\lambda}^{k+1}-\underline{\lambda}^k \rVert_{\infty} <
  10^{-3}.
\end{align*}
This condition is needed, due to underrelaxation, which could lead to a
non-changing active set, even though the solution is still changing. For
this illustrative example, simple underrelaxation is sufficient, but in
general, more sophisticated strategies are advised. To gain insight on
the condition number of $\widetilde{K}_h$ each nested conjugate gradient
iteration starts with a zero initial guess and stops at a relative residual
of $10^{-10}$. The experiment is conducted for varying mesh sizes
$h_t \simeq h_x$. The number of needed Newton as well as CG iterations are
given in Table \ref{tab:results}. Furthermore the table contains the ratio
of conjugate gradient iterations per Newton iterations, which stay at
reasonable values. In Figure \ref{fig:constrained_solution} we plot the
solution for different discretization parameters $n_x$ along the line
$x_1 = x_2 = x_3 = 0.51$ and compare it to the target function. A point
slightly off the center is chosen to exclude benevolent symmetry effects.
It is clearly visible, that the solution approaches the target function
until it reaches the upper bound. This is exactly the desired behavior.

\begin{table}[h!]
  \centering
  \caption{Numerical results for the constrained optimal control problem, with $n = n_t$.}\label{tab:results}
  \begin{filecontents*}{simulation_data/convergence_hist.csv}
Refinement, dof, n, NewtonIterations, TotalCG, relCG
0, 16, 2, 36, 36, 1
1, 256, 4, 36, 612, 17
2, 4096, 8, 36, 1296, 36
3, 65536, 16, 38, 2173, 57
4, 1.04858e+06, 32, 64, 3814, 60
\end{filecontents*}
\pgfplotstabletypeset[%
  col sep=comma, %
  columns={n, dof, NewtonIterations, TotalCG, relCG}, %
  every head row/.style={before row=\toprule,after row=\midrule}, %
    every last row/.style={after row=\bottomrule}, %
    columns/n/.style={column name={$\mathbf{n}$}, column type={r}}, %
    columns/dof/.style={column name=\textbf{DoF}, column type={r}, fixed}, %
  columns/NewtonIterations/.style={column name=\textbf{Newton iter.}, column type={r}}, %
  columns/TotalCG/.style={column name=\textbf{CG iter.}, column type={r}}, %
  columns/relCG/.style={column name=\textbf{CG/Newton}, column type={r}}, %
  ]{simulation_data/convergence_hist.csv}
\end{table}  

\paragraph{Acknowledgements} 
Part of this work has been supported by the Austrian
Science Fund (FWF) under the Grant Collaborative Research Center
TRR361/F90: CREATOR Computational Electric Machine Laboratory.

\end{document}